\documentclass{article}%
\usepackage{amsmath}%
\usepackage{amsfonts}%
\usepackage{amssymb}%
\usepackage{graphicx}%
\usepackage{color} 

\newtheorem{thm}{Theorem}

\newtheorem{cor}[thm]{Corollary}
\newtheorem{pro}[thm]{Proposition}

\def\Bbb{\mathbb}

\def\per{\mathop{\rm per}}

\def\wm<{\prec_{\tiny w}}
\def\wm>{\succ_{\tiny w}}
\def\lm<{\prec_{\mbox{\tiny $\log$}}}
\def\lm>{\succ_{\mbox{\tiny $\log$}}}
\def\wl<{\prec_{\mbox{\rm \tiny w$\log$}}}
\def\wl>{\succ_{\mbox{\rm \tiny w$\log$}}}

\def\psd<{\prec}
\def\psd>{\succ}
\def\psdgeq>{\succeq}

\def\n<{<}
\def\n>{>}

\def\proof{{\noindent \bf Proof. \hspace{.01in}}}
\newcommand{\qed}{\hspace{.1in} \vrule height 7pt width 5pt depth 0pt \medskip}

\begin{document}


\title{\Large On the number of vertices of the stochastic tensor polytope\thanks{To appear in Linear and Multilinear Algebra special issue dedicated to M. Marcus.}}

\author{Zhongshan Li${}^{\,\rm a}$
,\;
Fuzhen Zhang${}^{\,\rm b}$
,\;
Xiao-Dong Zhang${}^{\,\rm c}$
\\
\footnotesize{${}^{\,\rm a}$ Georgia State University, Atlanta, USA; zli@gsu.edu}\\
\footnotesize{${}^{\,\rm b}$ Nova Southeastern University, Fort Lauderdale, USA; zhang@nova.edu}\\
\footnotesize{${}^{\,\rm c}$ Shanghai Jiao Tong University, Shanghai, China;  xiaodong@sjtu.edu.cn}}


\date{}
\maketitle

\noindent In memory of Professor Marvin Marcus (1927-2016).

\bigskip
 \hrule
\bigskip

\noindent {\bf Abstract}
This paper is devoted to the study of lower and upper bounds for the number of vertices of the polytope of $n\times n\times n$ stochastic tensors (i.e., triply stochastic arrays of dimension $n$). By using  known results on   polytopes (i.e., the Upper and Lower Bound Theorems),
we present some new lower and upper bounds.
We show that the new upper bound is tighter than the  one recently
obtained by Chang, Paksoy and Zhang [Ann. Funct. Anal. 7 (2016), no.~3, 386--393] and also sharper than the one in Linial and Luria's [Discrete Comput. Geom. 51 (2014), no.~1, 161--170]. We demonstrate that the
analog of the lower bound obtained in such a way, however, is
no better than the existing ones.

\medskip


\bigskip
\hrule

\date{}
\maketitle


\section{Introduction}
Recall the well-known Birkhoff polytope $\omega_n$ of $n\times n$ doubly stochastic matrices: $\omega_n$ is the convex hull of all $n\times n$  permutation matrices.
As a polytope in $\Bbb R^{n^2}$, $\omega_n$ has dimension $(n-1)^2$ with $n^2$ facets, and $n!$ vertices.

Consider multi-arrays of higher dimension (here we focus on 3rd order). By a {\em stochastic tensor} (or {\em cube}) of dimension $n$  we mean a real $n\times n\times n$ hypermatrix (a.k.a. 3rd order tensor of dimension $n$)  $A=(a_{ijk})$ satisfying the conditions:
\begin{equation}
a_{ijk}\geq 0, \quad 1\leq i, j, k\leq n \label{eq10}
 \end{equation}
\begin{equation}
\sum_{i=1}^n a_{ijk}=1, \quad 1\leq  j, k\leq n \label{eq11}
 \end{equation}
\begin{equation}
\sum_{j=1}^n a_{ijk}=1, \quad 1\leq  i, k\leq n \label{eq12}
 \end{equation}
 \begin{equation}
 \sum_{k=1}^n a_{ijk}=1, \quad 1\leq  i, j\leq n\label{eq13}
 \end{equation}

Let $\Omega_n$ be the set of all $n\times n\times n$ stochastic tensors.
Consider each $A\in \Omega_n$ as an element of $\Bbb R^{n^3}.$
Then $\Omega_n$ is a subset of $\Bbb R^{n^3}$. Because it is an intersection
of a finite number of closed half spaces (as a linear equation $Ax=b$ is equivalent to
$Ax\leq b$ and $Ax\geq b$) and it is bounded, $\Omega_n$ is a polytope; that is,
it is generated by (or a convex hull of) finitely many points, i.e., extreme points or vertices of $\Omega_n$.
0-1 permutation tensors (i.e., stochastic tensors with 0 and 1 entries) are vertices of the polytope, but there are others in general (when $n>2$). For the case of $\Omega_3$, it is known that
there are 12 0-1 permutation tensors as vertices and 54 other vertices that are non  0-1
(see \cite[p.\,34]{Mayathesis03}). For a general $n$, the number of  0-1 permutation tensors of order 3 and dimension $n$ is the same as the number of $n\times n$ Latin squares (see \cite[p.~159]{Lint92} or \cite[Proposation 2.6]{CLN14}).

Note that the Birkhoff--von Neumann theorem for multistochastic tensors is investigated by Cui, Li and Ng \cite{CLN14}, Fischer and Swart \cite{FiSw85}, etc.,
and that a similar polytope of {polystochastic matrices} has been studied
by Gromova \cite{Gro74} (see also \cite[p.~64]{Bar02}), and by Brualdi and Csima \cite{BrCs75,BrCs75laa}. To our knowledge, Juflkat and  Ryser
\cite{JufRys68} are the first ones who directly addressed and studied the multidimensional matrices which now we call {\em tensors} (while general tensors have broader meanings depending contexts).

Although $\Omega_n$ is an analog of $\omega_n$ for higher order,
 $\Omega_n$ has many different geometric properties than $\omega_n$.
 The determination of the number and structures of the vertices of $\Omega_n$ is a very difficult problem. Estimation of the number of vertices  has been witnessed
 in three  ways:
{(1)}.  Combinatorial method via Latin squares.
Ahmed,  De Loera, and  Hemmecke (see \cite[Theorem~0.1]{Mayathesis03} or \cite[Theorem 2.0.10]{Mayathesis04})
 gave an explicit lower bound $\frac{(n!)^{2n}}{n^{n^2}}$.
{(2)}. Analytic (and direct) approach by using hyperplane and induction. Chang, Paksoy and Zhang \cite{ChangPZ16} recently showed an upper bound (see Theorem~\ref{Thm1} below). {(3)}. Computational geometry approach using known results on polytopes. Adopting this approach, we in this paper
present some upper and lower bounds of the number of vertices of $\Omega_n$ and compare the new bounds with the existing ones.

\begin{thm}\label{Thm1}
 Let $f_0(\Omega_n)$ be the number of vertices of the polytope $\Omega_n$.
Then
$$\frac{(n!)^{2n}}{n^{n^2}}\leq f_0(\Omega_n)\leq \frac{1}{n^3} \cdot {p(n)\choose n^3-1}, \;\; \mbox{where $p(n)=n^3+6n^2-6n+2$}$$
\end{thm}

The lower bound in Theorem \ref{Thm1} appeared in \cite{Mayathesis04} and \cite{Mayathesis03}; the upper bound  was recently obtained in \cite{ChangPZ16}. The
  bounds in Theorem \ref{Thm1} are very loose as we see in
Table~\ref{bounds}.
Note that the exact number of the vertices of $\Omega_4$ is unknown to our knowledge at this time.
For a general positive integer $n$, the
determination of the number
of vertices of $\Omega_n$ would be extremely difficult.

\begin{table}[h]
  \centering
  \begin{tabular}{c c c c} 
\hline
Case & lower & $f_0(\Omega_n)$ & upper \\ [0.5ex] 
\hline 
\\
$n=2$ & 1    & 2  & 21318  \\ [2ex]
$n=3$ & 2.37 & 66 &  $\frac{1}{27}{65 \choose 26}$
 \\ [2ex] 
 $n=4$ & 25.63 & ? &  $\frac{1}{64}{138 \choose 63}$
 \\[2ex]
\hline 
\end{tabular}
\caption{Lower and upper bounds in Theorem \ref{Thm1}}\label{bounds}
\end{table}

 \section{A sharper  upper bound}
 We use the standard terminology in convex polytope theory such as $d$-dimensional polytope and $i$-faces of a polytope. A 0-face is  a vertex (or an extreme point) and a $(d-1)$-face is called a facet; $f_i$ denotes the number of $i$-faces.
 One may refer to texts \cite{Bro83} and \cite{Zie95}  for the definitions of these terms. A fundamental question in the theory of convex polytope is the determination of the largest and the smallest numbers of $i$-faces, for instance,  the vertices,  of a polytope.

 Proved by McMullen in 1970 \cite{McM70}, the Upper Bound Theorem (UBT) is one of the most important results in the combinatorial theory of polytopes.  The UBT gives the maximum number of faces, say facets, of any $d$-polytope with a given number of vertices, and in a dual form, it gives the maximum number of vertices of any $d$-polytope with a given number of $k$-faces.  Further, the UBT establishes that the asserted maximum numbers are achieved by cyclic polytopes (see, e.g., \cite[p.\,16 and p.\,254]{Zie95}). As a consequence (see, e.g., \cite[p.\,90]{Bro83} and use duality),
 the number $f_0$ of vertices of a convex polytope of dimension $d$ with $f_{d-1}$ facets is bounded as follows:
\begin{equation}\label{upper}
f_0\leq  { {f_{d-1} - \lfloor \frac {d+1} {2} \rfloor } \choose { f_{d-1}-d}}  +  {{f_{d-1} - \lfloor \frac {d+2} {2} \rfloor } \choose { f_{d-1}-d} }
\end{equation}

Now we turn our attention to the number of vertices of the polytope  $\Omega_n$ of $n\times n \times n$ stochastic tensors. Regarded as a subset of $\Bbb R^{n^3}$,\;
 $\Omega_n$  defined  by  (\ref{eq10})-(\ref{eq13}) is  the same as the set of all
$x=(x_{ijk})\in \Bbb R^{n^3}$  satisfying
\begin{equation}
x_{ijk}\geq 0, \quad 1\leq i, \;j, \; k\leq n \label{eq0}
 \end{equation}
\begin{equation}
\sum_{i=1}^nx_{ijk}=1, \quad 1\leq j\leq n,\; 1\leq k\leq n \label{eq1}
 \end{equation}
 \begin{equation}
\sum_{j=1}^nx_{ijk}=1, \quad 1\leq i\leq n,\; 1\leq k\leq n-1 \label{eq2}
 \end{equation}
 \begin{equation}
\sum_{k=1}^nx_{ijk}=1, \quad 1\leq i\leq n-1,\; 1\leq j\leq n-1 \label{eq3}
\end{equation}

 Observe that  $n(2n-1) + (n-1)^2 = 3n^2 -3n+1$ independent line sum conditions are needed in defining $\Omega_n$ as a subset of $\mathbb R^{n^3}$. (In fact, the rank of the coefficient matrix of the linear equation system (\ref{eq1})-(\ref{eq3}) is $3n^2 -3n+1$.)
  Since each linearly independent linear equation reduces the dimension by 1, it follows that $\Omega_n$ as a polytope in $\Bbb R^{n^3}$ has dimension  $d=n^3- (3 n^2-3n +1) = (n-1)^3$. Alternatively, it can be computed as follows: Of the $n^3$ variables, $n^3-3n^2+3n-1=(n-1)^3$ variables are free (equivalently, we can take the variables in the $(n-1)\times (n-1)\times (n-1)$ cube in the lower-front-left corner as free variables).
\begin{figure}[h]
\begin{center}
\includegraphics[width=2in,totalheight=1.8in]{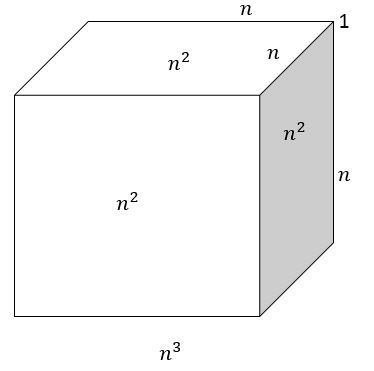}
\end{center}
\caption{Computing the dimension by using the cube}
\end{figure}

 As a polytope in $\Bbb R^{n^3}$,   $\Omega_n$ is of dimension $d=(n-1)^3$ and  has
 $f_{d-1}=n^3$ facets since each of the conditions in (\ref{eq0}) yields a facet.
 By Carath\'eodory's Theorem,
 any $n\times n$ doubly stochastic matrix is a convex combination of at most
 $(n-1)^2+1$ permutation matrices. Likewise, any $n\times n\times n$ 
 stochastic tensor is a convex combination of at most
 $(n-1)^3+1$ vertices (not necessarily permutation tensors).

  Applying the above result (\ref{upper}) to $\Omega_n$, with $d= (n-1)^3$ and $f_{d-1}=n^3$, we arrive at a tighter   bound for the number of vertices of $\Omega_n$.

 \begin{thm}\label{Thm2}
  Let $f_0(\Omega_n)$ be the number of vertices of the polytope $\Omega_n$.
Then
\begin{equation}\label{upper2}
f_0(\Omega_n)\leq
\left (\hspace{-.08in} \begin{array}{c}
n^3- \lfloor \frac{(n-1)^3+1}{2}\rfloor\\
 3n^2-3n+1
 \end{array} \hspace{-.08in} \right )+\left (\hspace{-.08in}  \begin{array}{c}
n^3- \lfloor \frac{(n-1)^3+2}{2}\rfloor\\
 3n^2-3n+1
 \end{array} \hspace{-.08in} \right  )
 \end{equation}
\end{thm}

\proof
The polytope  $\Omega_n$ is contained in  an $(n-1)^3$-dimensional affine subspace of the $n^3$-dimensional space $\Bbb R^{n^3}$ determined by the linear equality constraints that the sum of each row and of each column be one. Within this affine subspace, it is defined by $n^3$ linear inequalities, one for each coordinate of the tensor (hypermatrix), specifying that the coordinate be non-negative. Therefore, it has exactly $n^3$ facets
$F_{ijk}=\{x\in \Omega_{n} \mid x_{ijk}=0\}$, for $1\leq i, j, k\leq n$.  That is, $d=(n-1)^3$ and $f_{d-1}=n^3$ in (\ref{upper}).
This completes the proof.
\qed

We claim  that the upper bound in Theorem \ref{Thm2} is better (tighter) than the one in Theorem~\ref{Thm1}.
 That is, we show the following inequality.

\begin{pro}
 Let $n\geq 2$ be a positive integer. Then
\begin{equation*}\label{1}\left(\begin{array}{c}
n^3-\lfloor\frac{(n-1)^3+1}{2}\rfloor\\
3n^2-3n+1\end{array}\right)+\left(\begin{array}{c}
n^3-\lfloor\frac{(n-1)^3+2}{2}\rfloor\\
3n^2-3n+1\end{array}\right)< \frac{1}{n^3}
\left(\begin{array}{c}
n^3+6n^2-6n+2\\
n^3-1\end{array}\right)
\end{equation*}
\end{pro}

\proof
Direct computations for the cases $n=2, 3, 4, 5$  
show that the proposition holds. (In fact, for each case, the number on the right hand side is much larger than that on the left hand side.)  Now we assume $n\ge 6.$

 Let $a, b, k$ be positive integers such that $a\ge b$ and $k\ge 1$. Bear in mind the monotonicity properties that
 ${a\choose b} \leq {a+k\choose b}$ and that
 ${a\choose b} \leq {a\choose b+k}$ if $b+k\leq \lfloor \frac{a}{2} \rfloor$.

It follows that,  for $n\geq 6$ (which is used for the last inequality),
$$\left(\begin{array}{c}
n^3-\lfloor\frac{(n-1)^3+1}{2}\rfloor\\
3n^2-3n+1\end{array}\right)+\left(\begin{array}{c}
n^3-\lfloor\frac{(n-1)^3+2}{2}\rfloor\\
3n^2-3n+1\end{array}\right)$$
$$ <2\left(\begin{array}{c}
n^3\\
3n^2-3n+1\end{array}\right)<2\left(\begin{array}{c}
n^3\\
3n^2\end{array}\right)
$$
On the other hand, for $k\ge 1$,
\begin{equation}\label{2.1}\left(\begin{array}{c}
a+k\\
b+k\end{array}\right)
=\prod_{i=1}^k\frac{(a+i)}{(b+i)}\cdot \left(\begin{array}{c}
a\\
b\end{array}\right)>\left(\frac{a+k}{b+k}\right)^k\cdot \left(\begin{array}{c}
a\\
b\end{array}\right)\end{equation}
Hence (as $n\ge 6$), 
we have  $n^3+6n^2-6n+2\ge 2(6n^2-6n+3)$
and
\begin{eqnarray*}\label{3}
\lefteqn{\frac{1}{n^3}
\left(\begin{array}{c}
n^3+6n^2-6n+2\\
n^3-1\end{array}\right)}\\
& & =\frac{1}{n^3}
\left(\begin{array}{c}
n^3+6n^2-6n+2\\
6n^2-6n+3\end{array}\right)\\
& & > \frac{1}{n^3}
\left(\begin{array}{c}
n^3+6n^2-6n+2\\
3n^2+20\end{array}\right)
\ \ ({\rm by \ monotonicity}) 
\\
& & >
 \frac{1}{n^3}\left(\frac{n^3+6n^2-6n-18}{3n^2+20}\right)^{20} \cdot
\left(\begin{array}{c}
n^3+6n^2-6n-18\\
3n^2\end{array}\right)
\ \  ({\rm by \ (\ref{2.1})})  \\
 & & >\frac{1}{n^3}
\left(\frac{n}{3}\right)^{20} \cdot
\left(\begin{array}{c}
n^3\\
3n^2\end{array}\right)
\ \  ({\rm by \ monotonicity}) %
\\
& & >2
\left(\begin{array}{c}
n^3\\
3n^2\end{array}\right)
\ \
\end{eqnarray*}
This completes the proof. \qed

\begin{table}[h]
  \centering
  \begin{tabular}{c c c c} 
\hline
Case & $f_0(\Omega_n)$  & new up. bound   & old up. bound  \\ [0.5ex] 
\hline 
\\
$n=2$ &  2      &  2     &  21318  \\[2ex] 
$n=3$ &  66     &  10395  & $\frac{1}{27}{65 \choose 26}$  \\ [2ex] 
$n=4$ & ?     &  2 ${50 \choose 37}$  & $\frac{1}{64}{138 \choose 63}$ \\ [2ex] 
\hline 
\end{tabular}
  \caption{Comparison of upper bounds in Theorems 1 \& 2}\label{newupperbounds}
\end{table}

It is stated in \cite[p.~170]{LL14} that
$\Omega_n$ (i.e., $\Omega_n^{(2)}$ in \cite{LL14}) has fewer than $n^{3n^2}$  vertices. We show that our bound in Theorem \ref{Thm2} is tighter than this one also.

\begin{pro} Let $n\geq 2$ be a positive integer. Then
$$
\left(\begin{array}{c}
n^3-\lfloor \frac{(n-1)^3+1}{2}\rfloor\\
3n^2-3n+1\end{array}\right)+
\left(\begin{array}{c}
n^3-\lfloor \frac{(n-1)^3+2}{2}\rfloor\\
3n^2-3n+1\end{array}\right)<n^{3n^2}
$$
\end{pro}

\proof
Direct computations show that the assertion holds for $n=2, 3, \dots, 10$.
Let $n\ge 11$. Then $\frac{en^2}{3n^2-3n+1}<1$, where $e\approx 2.718$. Observe the fact that
 $$ \left(\begin{array}{c}
m\\
k\end{array}\right)<\left(\frac{me}{k}\right)^{k}\;\;
\mbox{for any positive integers $m$ and $k$,  $m\geq k$}
$$
(which follows from the Stirling's formula $\sqrt{2\pi} \,k^{k+\frac12}e^{-k}<k!$). We have
\begin{eqnarray*}
\lefteqn{\left(\begin{array}{c}
n^3-\lfloor \frac{(n-1)^3+1}{2}\rfloor\\
3n^2-3n+1\end{array}\right)+
\left(\begin{array}{c}
n^3-\lfloor \frac{(n-1)^3+2}{2}\rfloor\\
3n^2-3n+1\end{array}\right)}\\
& & <2\left(\begin{array}{c}
n^3\\
3n^2-3n+1\end{array}\right) \\
 & & <2\left(\frac{en^3}{3n^2-3n+1}\right)^{3n^2-3n+1}\\
 & & =2\left(\frac{en^2}{3n^2-3n+1}\right)^{3n^2-3n+1} \cdot n^{3n^2-3n+1}\\
& & < 2n^{3n^2-3n+1}<n^{3n^2} \quad  \qed
\end{eqnarray*}

Thus, we conclude that our new upper bound is tighter than the existing ones (to our knowledge). In fact, the left hand side is much smaller than the right hand side. For instance, for $n=10$, the quantity on the left is
$2\cdot {635\choose 271}\approx  9.8 \times 10^{186}$,
while on the right hand side, we have $10^{300}$.

\section{On the lower bound}

In the previous section, we saw that the upper bound for the number of vertices of $\Omega_n$ obtained by the McMullen Upper Bound Theorem (UBT) is {better}  than the known ones. A natural and analogous question to ask is: how about the lower bound? In this section we obtain a lower bound through a lower bound theorem and compare the lower bound derived this way to the existing ones.

Let $L_n$ denote the number of $n\times n$ Latin squares. $L_n$ can be computed by
$${\displaystyle L_{n}=n!\sum _{A\in B_{n}}^{}(-1)^{\sigma _{0}(A)}(_{\ \ n}^{\per A})}$$
where $B_n$ is the set of all 0-1 $n \times n$ matrices, $\sigma_0(A)$ is the number of zero entries in matrix $A$, and $\per A$ is the permanent of matrix $A$; see Shao and Wei \cite{ShaoWei92}.
Note that
$L_{n}\geq {\frac {\left(n!\right)^{2n}}{n^{n^{2}}}}$ (see, e.g., \cite[p.~162]{Lint92}).

 Since every Latin square can be interpreted as a 0-1 permutation tensor
(see, e.g., \cite[p.~159]{Lint92} or \cite[Proposition~2.6]{CLN14}) and every $n\times n\times n$ 0-1 permutation tensor is an extreme point of $\Omega_n$, we have
\begin{equation}\label{nLn}
 {\frac {\left(n!\right)^{2n}}{n^{n^{2}}}}\leq L_n \leq f_0(\Omega_n)
 \end{equation}

In contrast to the Upper Bound Theorem (UBT), the so-called Barnette Lower Bound Theorem (LBT) holds only for simplicial polytopes (see, e.g.,
\cite[p.~166]{Deza94b}).
Note that
the polytope of stochastic tensors is not simplicial in general. For example, when $n=3$, if $\Omega_3$ were simplicial, then according to Barnette's formula (see, e.g.,
\cite{Bar1973}) $f_{d-1}\geq (d-1)f_0-(d+1)(d-2)$, with $d=8$, $f_{d-1}=27$, $\Omega_3$ would have no more than 12 extreme points. This contradicts the fact that $\Omega_3$ has 66 vertices.
Nevertheless, a lower bound for the number of vertices of a general $d$-polytope with $f_{d-1}$ facets is provided
in \cite[Theorem~1.4]{Deza94a} (see also   \cite[Theorem 1.4]{Deza94b})
\begin{equation}\label{lower}
f_0\geq l_0^d(f_{d-1})
\end{equation}
where
\begin{equation}\label{eq:l}
l_0^d(x)=k \quad \mbox{if and only if} \quad
u_0^d (k-1)< x \leq u_0^d(k)\end{equation}
with $m=k-1$ and $k$ in
\begin{equation}\label{eq:u} u_0^d(m)=
{m-\lfloor \frac{d}{2}\rfloor -1\choose \lfloor \frac{d-1}{2}\rfloor }+
{m-\lfloor \frac{d-1}{2}\rfloor -1\choose \lfloor \frac{d}{2}\rfloor }
\end{equation}

For example, if $x=27$ and $d=8$, we want to find $k$ for which
$${k-6\choose 3}+{k-5\choose 4} <27\leq {k-5\choose 3}+{k-4\choose 4}$$
By a straightforward computation, we arrive at $k=11$. So, $l_0^8(27)=11$.

\begin{thm}\label{Thm4}
 Let $f_0(\Omega_n)$ be the number of vertices of the polytope $\Omega_n$.
Then
\begin{equation}\label{low}
  f_0(\Omega_n)\geq l_0^{(n-1)^3}(n^3)
 \end{equation}
 where
 $l_0^{(n-1)^3}(n^3)=k$ is such that (see (\ref{eq:l}) and (\ref{eq:u}))
 $$u^{(n-1)^3}_0(k-1)<n^3\leq u^{(n-1)^3}_0(k)$$
\end{thm}

\proof Applying
(\ref{lower}) with $d=(n-1)^3$ and $f_{d-1}=n^3$
  yields immediately the lower bound (\ref{low}).
\qed

In light of Theorems \ref{Thm1} and  \ref{Thm2} on the upper bound, we would naturally propose
\begin{equation*}
  l_0^{(n-1)^3}(n^3)\ge L_n \;\; (\mbox{which is $\ge  \frac{(n!)^{2n}}{n^{n^2}}$})
 \end{equation*}

It turns out that this is not the case in general; see the second part of Proposition \ref{pro8}.
To compare the lower bounds, we first  observe a fact that
if $a$ and $b$ are positive integers such that $a\ge b+2$ and $b\ge 2$, then
\begin{equation}\label{fact1}
\left(\begin{array}{c}
a\\
b\end{array}\right)\ge  \left(\begin{array}{c}
a\\
2\end{array}\right)
\end{equation}
This is because $\frac{a-b+r-2}{r}\ge 1$ when $a\geq b+2$ for $r=3, 4, \dots, b$,   and
\begin{eqnarray*}
\left(\begin{array}{c}
a\\
b\end{array}\right) & = & \frac{a(a-1)\cdots(a-b+1)}{1 \cdot 2\cdots b }\\
& = &  \left(\begin{array}{c}
a\\
2\end{array}\right)\cdot \frac{a-b+1}{3}\cdot\frac{a-b+2}{4}\cdots\frac{a-b+r-2}{r}\cdots \frac{a-2}{b}\\
 & \ge &  \left(\begin{array}{c}
a\\
2\end{array}\right)
\end{eqnarray*}

\begin{pro}
 Let $n\ge 4$. If  $k$ is a positive integer satisfying
$$
\left(\begin{array}{c}
(k-1)-\lfloor \frac{(n-1)^3}{2}\rfloor-1\\
\lfloor\frac{(n-1)^3-1}{2}\rfloor\end{array}\right)+
\left(\begin{array}{c}
(k-1)-\lfloor \frac{(n-1)^3-1}{2}\rfloor-1\\
\lfloor\frac{(n-1)^3}{2}\rfloor\;\;\end{array}\right)
$$
$$<n^3\leq $$
$$
\left(\begin{array}{c}
k-\lfloor \frac{(n-1)^3}{2}\rfloor-1\\
\lfloor\frac{(n-1)^3-1}{2}\rfloor\end{array}\right)+
\left(\begin{array}{c}
k-\lfloor \frac{(n-1)^3-1}{2}\rfloor-1\\
\lfloor\frac{(n-1)^3}{2}\rfloor\end{array}\right)
$$
then
$$ k=(n-1)^3+2$$
\end{pro}

\proof
{\bf Case 1}.  Suppose that  $n$ is even.
Then $$\left \lfloor\frac{(n-1)^3}{2}\right \rfloor
=\frac{(n-1)^3-1}{2},\quad  \left  \lfloor\frac{(n-1)^3-1}{2} \right \rfloor
=\frac{(n-1)^3-1}{2}$$
Let $p=\frac{(n-1)^3+3}{2}$.
Then $p\ge 15$ (for $n\geq 4$). We write
$$
\left(\begin{array}{c}
(k-1)-\lfloor \frac{(n-1)^3}{2}\rfloor-1\\
\lfloor\frac{(n-1)^3-1}{2}\rfloor\end{array}\right)+
\left(\begin{array}{c}
(k-1)-\lfloor \frac{(n-1)^3-1}{2}\rfloor-1\\
\lfloor\frac{(n-1)^3}{2}\rfloor\end{array}\right)
= 2\left(\begin{array}{c}
k-p\\
p-2\end{array}\right)$$
and
$$\left(\begin{array}{c}
k-\lfloor \frac{(n-1)^3}{2}\rfloor-1\\
\lfloor\frac{(n-1)^3-1}{2}\rfloor\end{array}\right)+
\left(\begin{array}{c}
k-\lfloor \frac{(n-1)^3-1}{2}\rfloor-1\\
\lfloor\frac{(n-1)^3}{2}\rfloor\end{array}\right)=2\left(\begin{array}{c}
k-p+1\\
p-2\end{array}\right)
$$
The inequalities in the proposition become
 \begin{equation}\label{Eq:eq1}
2\left(\begin{array}{c}
k-p\\
p-2\end{array}\right)<n^3\le 2\left(\begin{array}{c}
k-p+1\\
p-2\end{array}\right)\end{equation}
By (\ref{Eq:eq1}),  $k-p+1\ge p-2+1$, i.e., $k\ge 2p-2$.
If $k=2p-2,$  then (\ref{Eq:eq1}) becomes
$$2<n^3\le 2(p-1)=(n-1)^3+1$$ which is impossible for $n>1$.
Hence, $k\ge 2p-1.
$
We claim $k= 2p-1.
$
Suppose otherwise that $k\ge 2p$. Then $k-p\ge p-2+2$  and $p-2\ge 2$. By (\ref{fact1})
we have
\begin{eqnarray*}
2\left(\begin{array}{c}
k-p\\
p-2\end{array}\right)& \ge & 2\left(\begin{array}{c}
k-p\\
2\end{array}\right)=(k-p)(k-p-1)\ge (2p-p)(2p-p-1)\\
 & = & \frac{\big ((n-1)^3+3\big )\cdot \big((n-1)^3+1\big )}{4}>\frac{(n-1)^6}{4}\ge n^3
 \end{eqnarray*}
 contradicting (\ref{Eq:eq1}).
So,
$k=2p-1=(n-1)^3+2$. With such  $k$, (\ref{Eq:eq1}) becomes
$$
(n-1)^3+1<n^3 \leq
\left ( \frac{(n-1)^3+1}{2} \right )
\left ( \frac{(n-1)^3+3}{2} \right )$$
which holds for $n\geq 4$  by  direct verifications.
\medskip

\noindent
{\bf Case 2}. Suppose that $n\ge 5$ is odd.
Then
$$\left \lfloor\frac{(n-1)^3}{2}\right \rfloor=\frac{(n-1)^3}{2}, \quad  \left \lfloor\frac{(n-1)^3-1}{2}\right \rfloor=\frac{(n-1)^3-2}{2}$$
  Let $r=\frac{(n-1)^3+4}{2}=\frac{(n-1)^3}{2}+2\geq 5$.
  The inequalities in the proposition 
  becomes
 \begin{equation}\label{kr}
 \left(\begin{array}{c}
  k-r\\
  r-3\end{array}\right)+\left(\begin{array}{c}
  k-r+1\\
  r-2\end{array}\right)<n^3\le \left(\begin{array}{c}
  k-r+1\\
  r-3\end{array}\right)+\left(\begin{array}{c}
  k-r+2\\
  r-2\end{array}\right)
  \end{equation}
We claim that $k=2r-2=(n-1)^2+2$. To show this, we draw contradictions for both cases
(i) $k< 2r-2$ and (ii) $k> 2r-2$ .

 (i) If $k\le 2r-3$,  i.e.,  $k-r+1\le (r-3)+1,$ $ k-r+2\le (r-2)+1$, then
  \begin{eqnarray*}
  \left(\begin{array}{c}
  k-r+1\\
  r-3\end{array}\right)+\left(\begin{array}{c}
  k-r+2\\
  r-2\end{array}\right) & \leq   & (k-r+1)+(k-r+2)\\
    & = & 2k-2r+3\\
  & \le &  2(2r-3)-2r+3 \\
  & = & 2r-3
    =  (n-1)^3 +1<n^3
  \end{eqnarray*}
contradicting (\ref{kr}).

(ii) If  $k\ge 2r-1$, then $k-r\ge (r-3)+2, $ $ k-r+1\ge (r-2)+2 $.
  By (\ref{fact1}), we have
  \begin{eqnarray*}
  \left(\begin{array}{c}
  k-r\\
  r-3\end{array}\right)+\left(\begin{array}{c}
  k-r+1\\
  r-2\end{array}\right) & \ge & \left(\begin{array}{c}
  k-r\\
  2\end{array}\right)+\left(\begin{array}{c}
  k-r+1\\
  2\end{array}\right)\\
   & = & (k-r)^2\geq (r-1)^2 \\
   & > & \frac{(n-1)^6}{4}>n^3
   \end{eqnarray*}
  also contradicting (\ref{kr}).
 It follows that
  $k=2r-2=(n-1)^3+2.$ \qed

\begin{cor} Let $l_0^{(n-1)^3}(n^3)$ be defined by (\ref{eq:l}) and (\ref{eq:u}). Then
  $$l_0^{(n-1)^3}(n^3)=\left\{
  \begin{array}{ll}
  11, & for \ n=3\\
  (n-1)^3+2, & for \ n\ge 4\\
   \end{array}
  \right.$$
  \end{cor}


\begin{pro}\label{pro8}
Let $l_0^{(n-1)^3}(n^3)$ be defined by (\ref{eq:l}) and (\ref{eq:u}).

(1). If $n=3$ or $n=4,$ then
$$l_0^{(n-1)^3}(n^3)>\frac{(n!)^{2n}}{n^{n^2}}$$

(2). If $n\ge 5,$ then
$$l_0^{(n-1)^3}(n^3)< \frac{(n!)^{2n}}{n^{n^2}}$$
\end{pro}

\proof We compute as follows:\\

  $n=3: \quad \frac{(n!)^{2n}}{n^{n^2}}=\frac{64}{27}=2.37< l_0^{8}(27)=11 $\\

   $n=4: \quad \frac{(n!)^{2n}}{n^{n^2}}=\frac{6561}{256}=25.6< l_0^{27}(64)=29   $\\

      $n=5: \quad \frac{(n!)^{2n}}{n^{n^2}}=\Big (\frac{24^2}{5^3}\Big )^5=4.6^5> l_0^{64}(125)=66$\\

    $n=6: \quad  \frac{(n!)^{2n}}{n^{n^2}}=\Big (\frac{20}{6}\Big )^{12}>3^{12} >  l_0^{125}(216)=127$\\

    $n=7: \quad \frac{(n!)^{2n}}{n^{n^2}}=\frac{24^{14}30^{14}}{7^{35}}   =\Big (\frac{24}{7}\Big )^{14}\cdot \Big (\frac{30}{7}\Big )^{14}\cdot \frac{1}{7^7}>  l_0^{216}(343)=218$
\medskip

   For $n\ge 8$,  using the fact that $\sqrt{2n\pi}\left ( \frac{n}{e} \right )^n < n!$, 
   we have
   \begin{eqnarray*}
   (n-1)^3+2 & \le & (2n\pi)^n \\
     & \leq & (2n\pi)^n \cdot \left(\frac{n}{e^2}\right)^{n^2} \\
     & = & \frac{\left ( \sqrt{2n\pi}\, \big (\frac{n}{e} \big) ^n  \right )^{2n}}{n^{n^2}} \\
     & < &  \frac{(n!)^{2n}}{n^{n^2}}\quad \qed
     \end{eqnarray*}

We thus conclude that the lower bound for the number of the vertices of $\Omega_n$ through the lower bound theorem for a general polytope is looser than the existing ones when $n\geq 5$.

\section{Discussions}

We have seen that the lower and upper bounds for the number of vertices of $\Omega_n$
\begin{equation}\label{LU}
L_n \le f_0(\Omega_n)\leq U_n
\end{equation}
where $L_n$ is the number of $n\times n$ Latin squares and $U_n$ (see Theorem \ref{Thm2}) is
$$U_n=\left (\hspace{-.08in} \begin{array}{c}
n^3- \lfloor \frac{(n-1)^3+1}{2}\rfloor\\
 3n^2-3n+1
 \end{array} \hspace{-.08in} \right )+\left (\hspace{-.08in}  \begin{array}{c}
n^3- \lfloor \frac{(n-1)^3+2}{2}\rfloor\\
 3n^2-3n+1
 \end{array} \hspace{-.08in} \right  )$$

The bounds depending solely and explicitly on $n$ in (\ref{LU}) are the ``best" ones to our knowledge. However, as Table \ref{newupperbounds}
(for the case of upper bound) and numerical computations show that these bounds are still very loose. This is also seen in the proofs of the present results.
Better estimates have been called for (see \cite[Section~4]{LL14}).
Determination of the exact number of vertices (or even harder, $d$-faces) of a general polytope is an extremely difficult problem. The bounds we have obtained in fact are good for all polytopes of  dimension $(n-1)^3$ with
$n^3$ facets. A better estimate of the number of vertices specifically for the polytopes of stochastic arrays or Birkhoff type or $\Omega_n$ in particular  is desired.

 In their paper \cite{LL14} on the vertices of $d$-dimensional Birkhoff polytope, Linial and Luria use Latin squares and present a lower bound for the number of vertices of $\Omega_n$ (i.e., $\Omega_n^{(2)}$ in \cite[Theorem 1.5]{LL14}) with indeterminants $o(1)$:

\begin{thm}[Linial and Luria \cite{LL14}] The polytope $\Omega_n$ has at least $L_n^{\frac32-o(1)}$ vertices and $L_n=\Big (\big (1+o(1)\big )\frac{n}{e^2}\Big )^{n^2}.$ (Note: the two $o(1)$s may be different.)
\end{thm}

Thus, for the lower bound, we can write $L_n^{\frac32-o(1)}\leq f_0(\Omega_n)$. The shortcoming of this lower bound is that it contains
the indeterminants $o(1)$s; it gives an estimate but not the exact value for given $n$.
Note that
it is known \cite[p.~162, Theorem 17.3]{Lint92} that $L_n^{1/n^2}\sim e^{-2}n$ as $n\rightarrow \infty$.


\bigskip
\noindent
{\bf Acknowledgement.}
 Fuzhen Zhang is thankful to
Richard Brualdi for drawing reference \cite{LL14} to his attention.  The authors thank Chi-Kwong Li for his comments in the early stage of the project.  Fuzhen Zhang's work
 was partially supported by National Natural Science Foundation of China (NNSFC) No.~11571220 via Shanghai University. Xiao-Dong Zhang's work
  was partially supported by NNSFC No.~11531001, No.~11271256 and NNSFC-ISF Research Program No.~11561141001.


\end{document}